\documentclass[12pt]
{amsart}

\usepackage{amssymb}
\usepackage{amsmath}
\usepackage{array}

\newcommand{\Pf}{{\em Proof}. }
\newcommand{\EPf}{\hfill$\square$}

\newcommand{\Z}{\mbox{$\mathbf Z$}}
\newcommand{\R}{\mbox{$\mathbf R$}}
\newcommand{\C}{\mbox{$\mathbf C$}}

\newcommand{\SU}[1]{\mbox{$\mathbf{SU}(#1)$}}
\newcommand{\U}[1]{\mbox{$\mathbf{U}(#1)$}}
\newcommand{\SP}[1]{\mbox{$\mathbf{Sp}(#1)$}}
\newcommand{\SO}[1]{\mbox{$\mathbf{SO}(#1)$}}

\newcommand{\Spin}[1]{\mbox{$\mathbf{Spin}(#1)$}}
\newcommand{\G}{\mbox{$\mathbf{G}_2$}}

\newtheorem{thm}{Theorem}

\newtheorem{rmks}{Remarks}

\newtheorem{lem}{Lemma}

\newtheorem{conj}{Conjecture}

\title{Taut submanifolds}
\author{Claudio Gorodski}\thanks{Partially supported by
FAPESP and CNPq.}
\address{Instituto de Matem\'atica e Estat\'\i stica\\Universidade
de S\~ao Paulo\\Rua do Mat\~ao, 1010\\
S\~ao Paulo, SP 05508-090\\
Brasil}
\email{gorodski@ime.usp.br}

\date{\today}

\begin{document}

\begin{abstract}
This is a short, elementary survey article about 
taut submanifolds. In order to simplify the
exposition, we restrict to the case 
of compact smooth submanifolds of Euclidean or spherical 
spaces. Some new, partial results concerning taut $4$-manifolds
are discussed at the end of the text. 
\end{abstract}

\maketitle 

\section{Introduction}
This is a short, elementary survey article about 
taut submanifolds. In order to simplify the
exposition, we restrict to the case 
of compact smooth submanifolds of Euclidean or spherical 
spaces. Sections~2 through~4 collect basic definitions
and general results. In Section~5, we 
explain in greater detail some useful techniques
that we use in Section~6 to prove some new, partial
classification results about taut $4$-manifolds.

\section{The spherical two-piece property}

Before introducing the concept of tautness,
it is instructive to discuss the STPP. 
Let $M$ be a compact surface embedded in an Euclidean
sphere $S^m$. We say $M$ has the \emph{spherical two-piece 
property}, or STPP for short, if 
$M\cap B$ is connected whenever $B$ is a closed
ball in $S^m$~\cite{Ba}.
The STPP is equivalent to requiring that every Morse distance function
of the form 
\[ L_q:M\to\R,\qquad L_q(x)=d(x,q)^2 \]
has exactly one local minimum. Intuitively, this can be
seen by remarking that the sublevel sets of the $L_q$ 
are exactly the closed balls $B$. 
Moreover, by replacing $q$ by $-q$, one sees 
that these conditions imply that 
every Morse distance function $L_q$
also has exactly one local maximum. 
Let $b_i$ be the $i$th Betti number of $M$
where, unless explicitly stated, we will always 
use $\mathbf Z_2$ coefficients
to simplify the exposition, and let $\mu_i$ be the number 
of critical points of a Morse distance function $L_q$. 
From the Morse relations for the Euler characteristic,
\[ \chi = b_0-b_1+b_2 = \mu_0-\mu_1+\mu_2, \]
one finally deduces that the STPP is equivalent to the fact that
$\mu_i=b_i$ for all $i$ and for all $q$ such 
that $L_q$ is a Morse function. 

Since the STPP is defined in terms of intersections 
with closed balls, it is immediate that it
is a conformally invariant property of compact 
surfaces in $S^m$. Of course, one can also 
consider the STPP for compact surfaces of Euclidean 
space. Since stereographic projection is a conformal
transformation, the theories in $S^m$ 
and $\R^m$ are equivalent. 

It follows from results of Kuiper and Banchoff (see~\cite{Ba})
that a compact surface substantially embedded in $S^m$ with 
the STPP is a round sphere or a cyclide of Dupin in $S^3$,
or the Veronese embedding of a projective plane in $S^4$.
Recall that an embedding into a sphere is called
\emph{substantial} if its image does not lie in a sphere 
of smaller dimension. Of course, it is enough for us to consider
only substantial embeddings. 

We finally note that the STPP condition is 
also equivalent to requiring that
the induced homomorphism
\[ H_0(M\cap B,\mathbf Z_2)\to H_0(M,\mathbf Z_2) \]
in \u{C}ech homology is injective for every closed ball $B$. 
The use of \u Cech homology allows one to use 
all closed balls rather than only those
determined by level sets of distance functions that 
are Morse functions. 
If one prefers to use singular homology, than
the equivalent condition is reformulated in terms 
that the above homomorphism
in singular homology be injective for
\emph{almost} every closed ball. 

\section{Taut submanifolds}

Now we come to the main concept in this text.
Let $M$ be a compact submanifold embedded in $S^m$.
We say $M$ is \emph{taut} if
the induced homomorphism
\begin{equation}\label{homom}
H_i(M\cap B,\mathbf Z_2)\to H_i(M,\mathbf Z_2) 
\end{equation}
in \u{C}ech homology is injective for every closed ball $B$
and for all $i$.
This is equivalent to requiring that every Morse distance function
$L_q$ satisfies $\mu_i=b_i$ for all $i$, i.e.~$L_q$
is $\Z_2$-\emph{perfect}. Taut submanifolds were 
first considered by Carter and West~\cite{CW}. 
It is clear that tautness is equivalent to the STPP
in the case of surfaces. In the case of $3$-manifolds,
there is only a classification of the diffeomorphism
types of manifolds admitting taut embeddings, which was
given by Pinkall and Thorbergsson~\cite{P-Th}. They
remark that a geometrical classification, if possible, 
is expected to be a complicated problem since 
most of the examples already known depend on many parameters.
Their result is that there are seven 
types of such $3$-manifolds, namely 
$S^3$, $\R P^3$, $S^3/Q$,
$S^1\times S^2$, $S^1\times \R P^2$, $S^1 \times_h S^2$ and $T^3$.
Here $Q=\{\pm1,\pm i,\pm j,\pm k\}$ is the quaternion group, and
$h$ is a diffeomorphism of $S^1\times S^2$ that acts
on each factor as the antipodal map. A classification 
of $4$-manifolds admitting taut embeddings along similar 
lines seems feasible, and we will come back to this 
point in the last section. 

Next we present other classes of examples of taut
submanifolds. Clifford tori and the standard 
embeddings of projective spaces are important 
examples. We will not justify this assertion now, since 
these examples will be generalized below. 
It follows from the Chern-Lashof theorem~\cite{C-L1}
that the sphere $S^n$ can only be tautly embedded in Euclidean
space as a round hypersphere in an affine subspace. 
This result also admits an alternative proof for $n\geq2$
by noting that the homology of $S^n$ is trivial
except in dimensions $0$ and $n$, which implies
that a Morse distance function to a taut
embedding of $S^n$ can have 
only critical points of index $0$ and $n$,
and so a focal point is necessarily of
multiplicity $n$, which yields that $S^n$ is 
umbilic according to an argument sketched below
in Section~5 (compare~\cite{No-Ro}).
Moving towards manifolds with more complicated topology,
Cecil and Ryan proved that a taut $n$-dimensional compact  
hypersurface of $S^{n+1}$ with the
same homology as $S^k\times S^{n-k}$ has precisely two principal
curvatures at each point, and the principal curvatures are constant
along the corresponding curvature distributions~\cite{C-R2}.
They called such a hypersurface a (high-dimensional) 
\emph{cyclide of Dupin}. 
Thorbergsson considered the case of a taut compact submanifold
of dimension $2k$ that is $(k-1)$-connected but not 
$k$-connected, and showed that it must be 
either a cyclide of Dupin diffeomorphic to 
$S^k\times S^k$ or the standard 
embedding of a projective plane over one of the
four normed division algebras~\cite{Th7}. 

In another development, Bott and Samelson proved 
that the orbits of the isotropy representations 
of the symmetric spaces, sometimes called
\emph{$s$-representations}, are tautly embedded
by explicitly constructing cycles forming a basis in $\mathbf{Z}_2$-homology
for these orbits~\cite{B-S}.
These orbits are known as the \emph{generalized
flag manifolds}. The generalized flag manifolds 
are homogeneous examples of another very important,
more general class of submanifolds, called
\emph{isoparametric submanifolds}. Hsiang, Palais 
and Terng proved that isoparametric submanifolds and their 
focal submanifolds are taut~\cite{H-P-T}. 

\section{Taut representations}

Most of the known examples of taut embeddings are 
homogeneous spaces. On the other hand, Thorbergsson derived some necessary 
topological conditions for a homogeneous space 
to admit a taut embedding, which allowed him to exhibit
some examples of homogeneous spaces which cannot be
tautly embedded, e.g.~the lens spaces distinct from
the real projective space~\cite{Th4}. 

Thorbergsson and I approached the 
problem of classifying taut compact submanifolds 
of Euclidean space which are extrinsically 
Riemannian homogeneous by first studying \emph{taut
representations}, namely those representations
of compact Lie groups 
all of whose orbits are tautly embedded. 
The only examples known at that time were 
the $s$-representations, and many 
proofs had been given of the tautness
of special cases of generalized flag manifolds where the arguments were 
easier. In~\cite{G-Th} (see also~\cite{G-Th0}), we classified the taut
irreducible representations of the compact Lie groups.
It turns out that the classification
includes three families of
representations that are not $s$-representations, thereby
supplying many new examples of 
tautly embedded homogeneous spaces. These families
are given by the following table, where $n\geq2$.
\[ \begin{array}{|c|c|}
\hline
\SO2\times\Spin9 & \mbox{(standard)}\otimes_{\mathbf R}\mbox{(spin)} \\
\U2\times\SP n & \mbox{(standard)}\otimes_{\mathbf C}\mbox{(standard)} \\
\SU2\times\SP n & \mbox{(standard)}^3\otimes_{\mathbf H}\mbox{(standard)} \\
\hline
\end{array}\]
It is worth noting that these representations
are exactly those representations of cohomogeneity 
three that are not $s$-representations. 
In~\cite{G-Th2}, we showed that the orbits of these
representations also admit cycles of Bott-Samelson type.

The proof of the classification theorem in~\cite{G-Th}
is long and intrincate. It starts with a remark by Kuiper
that implies that a taut irreducible representation has
the property that the second osculating spaces of
all of its nontrivial orbits coincide with the ambient space~\cite{Ku2}. 
We call representations with this property,
irreducible or not,  
\emph{of class $\mathcal O^2$}. The class
$\mathcal O^2$ is much more easy to deal with since
it involves an infinitesimal condition. Indeed, 
we establish necessary upper bounds for the Dadok invariant 
$k(\lambda)$ (which is an integral algebraic invariant
of an irreducible representation,
see~\cite{D}) for irreducible representations of class $\mathcal O^2$,
which allows us to reduce a lot the size of that class
in such a way that the remaining cases are treated 
with geometrical methods that we develop in the
second part of the proof. 

There are three main strategies involved in these
methods. The first one often works a kind of induction.
It states that any slice representation of a taut
representations inherits the property of being taut. 
(Recall that the \emph{slice representation} of a 
representation at a point $p$ 
is the representation induced by the isotropy subgroup at $p$  
on the normal space to the orbit through $p$ at $p$.)
The second one might be called ``the fundamental
result about taut sums''. It relates the topology 
of the orbits of a taut reducible representation to 
the topology of the orbits of its summands. Finally,
the third strategy invokes a reduction principle 
in transformation groups which \emph{grosso modo}
reduces the task of deciding whether a given representation
is taut or not to the study of the tautness of a much simpler
representation with trivial principal isotropy subgroup. 

The combination of the above three strategies is also
effective in the classification of taut reducible
representations of compact \emph{simple} Lie groups
which I completed in~\cite{Gor2}. In this case,
%
%
since every orbit of a summand of a taut reducible
representation is also an orbit 
of the sum, it follows from the classification
in the irreducible case that 
every irreducible summand of a taut reducible representation
is an $s$-representation, and we need only 
to decide which sums of $s$-representations remain taut.
Each compact simple Lie 
group admits few $s$-representations, so that 
the analysis via the geometrical methods
discussed in the previous paragraph  
can be effectively carried out and the final result
is given in the following table. 
{\small
\[ \begin{array}{|c|c|c|}
\hline
\SU n,\;n\geq3 & \C^n\oplus\cdots\oplus\C^n & \mbox{$k$ copies, where $1<k<n$} \\
\hline
\SO n,\;n\geq3, n\neq4 & \R^n\oplus\cdots\oplus\R^n & \mbox{\textrm{$k$ copies, where $1<k$}} \\
\hline
\SP n,\;n\geq1 & \C^{2n}\oplus\cdots\oplus\C^{2n} & \mbox{\textrm{$k$ copies, where $1<k$}} \\ 
\hline
\G & \R^7\oplus\R^7 & \mbox{---} \\
\hline
\Spin6 & \R^6\oplus\C^4 & \R^6=\mathrm{(vector)},\;\C^4=\mathrm{(spin)} \\
\hline
       & \R^7\oplus\R^8 & \\
\Spin7 & \R^8\oplus\R^8 & \R^7=\mathrm{(vector)},\;\R^8=\mathrm{(spin)} \\
       & \R^8\oplus\R^8\oplus\R^8 & \\
       & \R^7\oplus\R^7\oplus\R^8 & \\
\hline
       & \R^8_0\oplus\R^8_+ & \\
\Spin8 & \R^8_0\oplus\R^8_0\oplus\R^8_+ & \R^8_0=\mathrm{(vector)},\;
\R^8_+=\mathrm{(halfspin)} \\
       & \R^8_0\oplus\R^8_0\oplus\R^8_0\oplus\R^8_+ & \\
\hline
\Spin9 & \R^{16}\oplus\R^{16} & \R^{16}=\mathrm{(spin)} \\ 
\hline
\end{array}\]
}

We make two simple remarks regarding the representations 
appearing in this table. The first one 
is that in the case of the first representation in the table,
any number of summands $\C^n$ in the sum can be replaced 
by the dual representation $(\C^n)^*$,
and the resulting representation remains taut. 
The second one is that 
the representations of $\Spin8$ are listed up to composition 
with an outer automorphism of the Lie group, so
the pair $(\R^8_0,\R^8_+)$ appearing in the list can be replaced by any 
pair of inequivalent $8$-dimensional
representations of $\Spin8$, and the resulting 
representations for $\Spin8$ will still
be taut.

Despite these results, the classification of taut 
extrinsically Riemannian homogeneous submanifolds of Euclidean space 
is still far from complete. For example, we still do not know
what the taut reducible representations of the nonsimple
groups are. 

\section{The Morse index theorem, proper 
Dupin submanifolds, and Ozawa's theorem}

In this section, we review a collection of
methods related to taut submanifolds 
that were used in~\cite{P-Th} and will  
be useful here. 

Let $M$ be a compact 
smooth $n$-manifold substantially embedded
in $S^{n+k}$. Later we will specify to the case in which 
$M$ is tautly embedded.
Denote by $N^1(M)$ the unit normal bundle of $M$
in $S^{n+k}$. 
A focal point of $M$ in $S^{n+k}$ is a critical value
of the restriction of the exponential map of $S^{n+k}$
to $N^1(M)$. According to this definition, 
$-p$ is a focal point of $M$ for every $p\in M$
if $k>1$. 
This kind of focal point is uninteresting for us
since it does not relate to the geometry of $M$,
but rather to the fact that any pair of antipodal points 
are conjugate in the sphere. Therefore in the following
we will completely disconsider this kind of focal point.
All the other focal points of $M$ can be described as follows.       
For $\xi\in N^1(M)$,
the set of focal points of $M$ lying in the normal ray defined by $\xi$ 
is discrete and in correspondence with
the principal curvatures of $M$ 
relative to the Weingarten operator $A_\xi$
in such a way that for a focal distance $d\in(0,\pi)$ there corresponds 
an eigenvalue $\cot d$ of $A_\xi$.

Now it can be checked that a point $q\in S^{n+k}$ is a focal point of $M$ 
if and only if the distance function $L_q$ is not a Morse function. 
In this case, $q$ is a focal point relative to 
some $\xi\in N^1_p(M)$, $p$ is a degenerate critical point of 
$L_q$, and the multiplicity of $q$ as a 
focal point relative to $\xi$
equals the nullity of the Hessian 
of $L_q$ at $p$, which is also the the multiplicity 
of the corresponding principal curvature of $A_\xi$.
For a nonfocal point $q$, 
$L_q$ is a Morse function and the Morse index theorem asserts 
that the index of a critical point $p$  
of $L_q$ equals the sum of the multiplicities of the 
focal points of $M$ lying in the normal geodesic segment $\overline{pq}$.    

It follows that the sum of the multiplicities of the 
focal points lying in any open half great circle joining $p$ and $-p$  
is $n$. Since $A_{-\xi}=-A_{\xi}$ and $-\cot d = \cot(\pi-d)$,
we also have that for each focal distance $d\in(0,\pi)$ 
in the direction of $\xi$, there is a corresponding 
focal distance $\pi-d$ in the direction of $-\xi$ 
with the same multiplicity. Let $n(\xi)$ denote the number
of distinct focal points lying in the open half great circle 
specified by $\xi$.  
We have that $n(\xi)$, as a function on $N^1(M)$, is lower
semicontinuous and   
there is an open and dense subset $\Omega$ of $N^1(M)$
where $n(\xi)$ is locally constant and maximal.
A normal vector $\xi\in N^1(M)$ is called
\emph{regular} if it belongs to $\Omega$, 
and \emph{singular} otherwise. 

We suppose henceforth that $M$ is tautly embedded
in $S^{n+k}$. 
 
Now the condition that every unit normal vector is regular
is equivalent to the condition that $M$ be
a \emph{proper Dupin} submanifold of $S^{n+k}$~\cite{Pi1,Th3}. 
The latter means that the principal curvatures are 
constant along the corresponding curvature 
surfaces and the number of distinct principal curvatures
is constant. 
In this case, there is a well defined 
ordered sequence $(m_1,\ldots,m_g)$ of multiplicities 
of focal points of $M$ along the normal ray defined by
an arbitrary $\xi$.
Here $g$ is the (constant) number of distinct principal curvatures 
of $M$. 

If $M$ is taut and there exists a unit normal vector $\xi\in N^1(M)$
such that $n(\xi)=1$, then $M$ is umbilic, and hence, 
a great hypersphere $S^n$ in $S^{n+1}$. This can be seen
as follows. Let $f$ be the first focal point of $M$ in the
direction of $\xi$. The multiplicity of $f$ is $n$
by hypothesis. Consider the geodesic segment $\overline{pf}$,
the open half great circle $\gamma$ containing it, and a point
$q$ belonging to $\gamma$. 
The Morse index theorem implies that 
the index of $p$ as a critical point of $L_q$ is $0$ 
or $n$, according to whether 
$q$ occurs between $p$ and $f$ or past $f$. It follows that
$p$ is a local minimum or a local maximum of $L_q$, respectively.
By tautness of $M$, there can be only one local minimum 
(resp.~local maximum), so this is also a global minimum 
(resp.~global maximum).
Letting $q$ tend to $f$ in both cases shows that $M$ is contained 
in the hypersphere of center $p$ and radius $d(p,f)$,
and hence $M$ coincides with that hypersphere. 

As was said above, if $M$ is taut and $q$ is a focal point, then 
$L_q$ is not a Morse function, but there is a very useful theorem 
by Ozawa~\cite{Oz} that asserts that 
$L_q$ is a $\Z_2$-perfect Morse-Bott function and the connected
components of the critical set of $L_q$ are compact smooth
taut submanifolds of $S^{n+k}$. Recall 
that a function $F$ is called a \emph{Morse-Bott}
function if the connected components of its 
critical set are smooth submanifolds along each of which 
the nullity of the Hessian of $F$ is constant 
and equal to the dimension of the component~\cite{Bo2}. 
In this case, the index of $F$ along a critical manifold 
is by definition the index of the Hessian 
of $F$ restricted to the normal bundle of that
critical manifold.
A Morse-Bott function $F$ on $M$ is \emph{perfect}
if it satisfies the 
Morse-Bott equalities, namely 
the Poincar\'e polynomials of $M$ and those of the
critical manifolds $C_i$ of $F$ are related by the formula
\[ \mathcal P_t(M) = \sum_i \mathcal P_t(C_i)\,t^{\mathrm{ind}C_i}. \]
If $M$ is taut, it follows from Ozawa's theorem and the 
Morse-Bott equalities that the $\Z_2$-homology of the critical manifold 
corresponding to the set of minima (resp.~maxima) 
of a distance function $L_q$ injects
into the $\Z_2$-homology of $M$. Such critical sets are called
\emph{top sets} or \emph{top cycles}, or yet \emph{top circles},
\emph{top tori}, etc., when their topology is specified, 
following terminology  
introduced by Kuiper and expanded by Pinkall and Thorbergsson. 
For each $\xi\in N^1(M)$, we denote by $T(\xi)$ the top
cycle which is the set of minima of $L_f$, where
$f$ is the first focal point in the normal ray
defined by $\xi$.

Another interesting consequence of Ozawa's
theorem is that the definition of a taut submanifold
can be restated so as to require that the 
homomorphism~(\ref{homom}) in \emph{singular} homology be injective for
\emph{all} closed balls $B$.

\section{Taut $4$-manifolds}

In this section we present some new, partial
results about the classification of compact 
smooth manifolds of dimension four that admit 
taut embeddings. 
We start with the two following results,
the first one of which has already been noticed by
Thorbergsson. 

\begin{thm}\label{b_1=0}
A compact four-dimensional smooth taut 
submanifold $M$ with vanishing first Betti number 
is diffeomorphic to $S^4$, $S^2\times S^2$ or $\mathbf CP^2$.
\end{thm}

\Pf If $b_2=0$, then $M$ is a $\Z_2$-homology sphere,
and we have already seen that it must be a hypersphere
$S^4$ in $S^5$. 

Next, we remark that a taut submanifold 
with $b_1=0$ must be simply connected. 
This is because it admits a $\Z_2$-perfect Morse
function, and such a function has no critical points 
of index one (compare Lemma~4.11(1) in~\cite{C-E}). 

If $b_2\neq0$, it follows from the preceeding remark 
that $M$ is a taut $4$-manifold that is $1$-connected
but not $2$-connected. Therefore Theorem~A in~\cite{Th7} says that either
$M$ is a cyclide of Dupin of type $S^2\times S^2$ in $S^5$, or
$M$ is diffeomorphic to $\C P^2$ and sits inside $S^7$. 
Recall that the standard embedding of $\C P^2$  
in $S^7$ is taut. \EPf

\begin{thm}\label{b_2=0}
A compact four-dimensional smooth taut 
submanifold $M$ with vanishing second Betti number 
is diffeomorphic to $S^4$ or $S^1\times S^3$.
\end{thm}

\Pf By Theorem~\ref{b_1=0}, we may assume that $b_1\neq0$. 
Then, by hypothesis and $\Z_2$-Poincar\'e duality, exactly four 
$\Z_2$-homology groups of $M$ are nonzero. It follows from 
Theorem~3.2 in~\cite{Oz} that $b_1=1$ (this also
follows from Theorem~11 in~\cite{He}). 

We will use some basic facts about the topology 
of $4$-manifolds~\cite{G-S}. 
Since there is a Morse function on $M$ with one critical point 
of each index $0$, $1$, $3$, $4$, there is a
handle decomposition 
$M=h_0\sqcup h_1\sqcup h_3\sqcup h_4$, where each $h_i$ an $i$-handle
attached via an attaching map. 

If $M$ is orientable, then $h_1$ is an orientable handle,
$h_0\sqcup h_1$ and $h_3\sqcup h_4$ are both 
diffeomorphic to $S^1\times D^3$, and $M$ is 
obtained by glueing two copies of $S^1\times D^3$ along their
boundaries by a diffeomorphism $f$ of $S^1\times S^2$. 
The diffeomorphism class of $M$ only depends on the 
diffeotopy class of $f$. It has been shown that the
diffeotopy group of $S^1\times S^2$ is isomorphic to 
$\Z_2\oplus\Z_2\oplus\Z_2$, 
where each generator extends to a diffeomorphism of $S^1\times D^3$
\cite{Gl}. It follows easily that there is at most 
one manifold in this case, which plainly is $S^1\times S^3$. 

If $M$ is nonorientable, then $h_1$ is a nonorientable
handle, $h_0\sqcup h_1$ and $h_3\sqcup h_4$ are both 
diffeomorphic to $S^1\tilde\times D^3$, the twisted 
$3$-disk bundle over $S^1$, and $M$ is 
obtained by glueing two copies of $S^1\tilde\times D^3$
along their boundaries
via a diffeomorphism of $S^1\tilde\times S^2$.
Every diffeomorphism of $S^1\tilde\times S^2$
is diffeotopic to one that extends to $S^1\tilde\times D^3$
\cite{K-R}. It follows as above that $M$ must be 
diffeomorphic to $S^1\tilde\times S^3$, but the 
following lemma shows that this case cannot occur. 
This finishes the proof. \EPf

\begin{lem}
We have that $S^1\tilde\times S^3$ does not admit a taut
embedding.
\end{lem}

\Pf Suppose, on the contrary, that $M=S^1\tilde\times S^3$ 
admits a taut and substantial embedding into a sphere $S^{4+k}$.
Since $b_2=0$, any top set must have vanishing second 
Betti number. It follows that
$M$ cannot admit a $2$-dimensional top set,
and any $3$-dimensional top set has the $\Z_2$-homology of $S^3$. 
Therefore, any top set of $M$, being taut, is diffeomorphic to $S^1$ 
or $S^3$. 
This implies that there is no $\xi\in N^1(M)$ 
with $n(\xi)=4$. For otherwise, 
if $q$ is a point lying in between the second and third 
focal points in the normal ray defined by $\xi$, 
and $p$ is the foot point of $\xi$, then $p$
is a critical point of index $2$ of $L_q$, but 
$b_2=0$. 

Now, for every $\xi$,  $n(\xi)=2$ or $3$. 
We claim that $M$ is proper Dupin. If not, 
there exist a singular $\bar\xi$ with $n(\bar\xi)=2$
and a sequence $\xi_i$ converging to $\bar\xi$
with $n(\xi_i)=3$. By replacing $\bar\xi$ and the $\xi_i$
by their opposites if necessary, we may assume that
the top set $T(\bar\xi)$ is diffeomorphic to $S^3$. 
The $T(\xi_i)$ are round circles representing 
nontrivial one-dimensional homology classes of $M$. A
subsequence of the $T(\xi_i)$ converges to a circle
in $T(\bar\xi)$ that is homologically nontrivial, and this is
a contradiction (compare the proof of
Theorem~2.2 in~\cite{Oz}). 

Now $M$ is proper Dupin. Suppose first that $n(\xi)=2$
for all $\xi$. We cannot have $k=1$, since $M$ is not
a cyclide of Dupin, so $k>1$. Note that the multiplicities
of the principal curvatures are $1$ and $3$.
Let $T$ be a small tube around
$M$. Then $T$ is a proper Dupin hypersurface with 
three distinct principal curvatures of multiplicities
$1$, $3$ and $k-1$. 
It is known that a proper Dupin
hypersurface with $3$ distinct principal curvatures
has all multiplicities equal~\cite{Mi}, so this case
cannot occur. 

Suppose now that $n(\xi)=3$ for all $\xi$. 
A small tube $T$ around $M$ is a 
proper Dupin hypersurface with  
$4$ distinct principal curvatures of multiplicities
$(1,2,1,k-1)$. It is known that in the case of 
$4$ distinct principal curvatures there are at most 
two different multiplicities~\cite{G-H}. It follows that $k=3$.
Let $B_+$ and $B_-$ be the two components of 
$S^7\setminus T$ and assume that $M\subset B_+$. It follows
from the Mayer-Vietoris sequence that
\[ H_3(T;\Z) = H_3(B_+;\Z)\oplus H_3(B_-;\Z). \]
Since $B_+$ is homotopy equivalent to $M$, we see
that $H_3(T;\Z)=H_3(M;\Z)\oplus H_3(B_-;\Z)=\Z_2\oplus H_3(B_-;\Z)$
(see~\cite{Hat}, p.~238). On the other hand, by~\cite{G-H},
$H_3(T;\Z)=\Z\oplus\Z$. This is a contradiction. Thus 
there can be no taut embedding. \EPf

\begin{rmks}
\rm (i) As a consequence of Theorem~\ref{b_2=0},
a compact smooth $4$-manifold with the same 
homology as $S^1\times S^3$ which is tautly embedded
in arbitrary codimension is diffeomorphic 
to $S^1\times S^3$ (compare Theorem~2 in~\cite{C-R2}
and the next item (ii)). 

(ii) According to Theorem~3.1 in~\cite{Oz}, the codimension
of a taut and substantial embedding of $S^1\times S^3$
into a sphere is either~$1$ or~$3$. Both codimensions can be realized.
In fact, $S^1\times S^3$ can be realized as a cyclide
of Dupin in codimension one, and it is shown in~\cite{C-R}
that any taut embedding of $S^1\times S^3$ with codimension one is 
such a cyclide, and they are all M\"obius equivalent to 
a tube around a circle in $S^5$. Moreover, 
$S^1\times S^3$ can also be tautly embedded in $S^7$
as a singular orbit of the isotropy representation of
the Grassmann manifold $\SO6/(\SO2\times\SO4)$. 

(iii) We have that $S^1\times_h S^3$, where $h$ acts
by the antipodal map on each factor, is diffeomorphic to
$S^1\times S^3$. In fact, $S^1\times_h S^3\approx \U1\times_{\mathbf Z_2}
\SU2=\U2$ and $\det:\U2\to\U1$ is a principal $\SU2$-bundle
with a global section $\sigma:\U1\to\U2$ given by 
\[ \sigma(e^{i\theta})
=\left(\begin{array}{cc} e^{i\theta}&0\\0&1\end{array}\right). \] 

(iv) According to Theorem~8 in~\cite{He}, $S^2\tilde\times S^2\approx
\C P^2\#\overline{\C P^2}$ cannot be tautly embedded into a sphere.
\end{rmks}

Regarding the following conjecture,
note that the principal orbits of the isotropy representation
of the Grassmann manifold $\SO5/(\SO2\times\SO3)$
are proper Dupin hypersurfaces in $S^5$
diffeomorphic to $S^1\times\mathbf RP^3$.

\begin{conj}\label{proper-dupin}
A compact embedded proper Dupin hypersurface in $S^5$
with four distinct principal curvatures is diffeomorphic
to $S^1\times\mathbf RP^3$.
\end{conj}

We next present a proof of this conjecture that uses
that the Poincar\'e conjecture in dimension~$3$ is true
(there has been recent progress by Perelman
on this conjecture, see~\cite{Mo}). 
Let $M$ be a compact embedded 
proper Dupin hypersurface in $S^5$
with four distinct principal curvatures. 
Then the multiplicities of focal points
along a normal ray are $(1,1,1,1)$. 
For a fixed unit normal vector field, 
let $F_+$ (resp.~$F_-$) denote the set of focal points
corresponding to the largest (resp.~smallest)
principal curvature, and let $f_{\pm}:M\to F_{\pm}$ be 
the focal maps. Then these are circle bundles
whose fibers consist of the circles of curvature.
It follows from the proof of Proposition~3.5
in~\cite{G-H} that $\pi_1(F_+)=\Z_2$ and $\pi_1(F_-)=\Z$
or $\pi_1(F_+)=\Z$ and $\pi_1(F_-)=\Z_2$, and we may assume 
we are in the first case just by replacing the normal field 
by its opposite if necessary. Now the universal covering 
manifold of $F_+$ is a simply connected compact 
$3$-manifold, hence diffeomorphic to $S^3$ 
if the Poincar\'e conjecture is true. 
The quotient of $S^3$ by a fixed point free
involution is diffeomorphic to $\R P^3$ according
to~\cite{Liv}, so $F_+$ is diffeomorphic 
to $\R P^3$. So far we know that $M$ is a circle
bundle over $\R P^3$. 
Besides the trivial one, there are exactly two other
circle bundles over $\R P^3$, one oriented, one nonoriented
(in fact, the associated $2$-plane bundles are $2L$ and $L\oplus T$,
where $T$ is the trivial line bundle and $L$ is the
tautological line bundle~\cite{Le}). 
Let $E_1$, $E_2$ respectively denote the total spaces of these 
bundles. Using the Gysin sequence
in the oriented case, one easily computes that
$H^2(E_1;\Z)=0$, and hence, $H_2(E_1;\Z)=0$. 
Using the Serre homology
spectral sequence in the nonoriented case, 
it is also not difficult to see that 
$H_1(E_2;\Z)$ is a $\Z_2$-extension 
of $\Z_2$. In any case, $E_1$ and $E_2$
do not have the integral homology of the trivial 
circle bundle $\R P ^3\times S^1$.
But $M$ must have the integral homology
of $\R P ^3\times S^1$ according to~\cite{G-H}. It follows that 
$M$ is diffeomorphic to $\R P ^3\times S^1$. 

\medskip

It is easy to construct other examples of 
taut embeddings of compact 
smooth $4$-manifolds
that are not covered by the previous results,
so there is still a lot of work 
to be done in order to achieve a complete 
classification.
In the case of taut hypersurfaces in $S^5$
with no more than three distinct principal
curvatures at each point, 
we have the following theorem~\cite{Ni}. 

\begin{thm}[Niebergall]
Let $M$ be a compact smooth
taut hypersurface in $S^5$. If $M$ has at most 
three distinct principal curvatures
at each point, then $M$ is diffeomorphic 
to one of the following $4$-manifolds:
$S^4$, $S^3\times S^1$, $S^2\times S^2$,
$T^2\times S^2$, or the nontrivial 
$S^2$-bundle over $\R P^2$.
\end{thm}

Niebergall proves this theorem by
first remarking that there is an open and dense
subset $M_0$ of $M$ where the multiplicities
of the principal curvatures are locally constant. 
Then each connected component of $M_0$
is a piece of a proper Dupin hypersurface 
with three distinct principal curvatures in~$S^5$. 
By the main result of~\cite{Ni}, each
connected component of $M_0$ is then a \emph{reducible}
proper Dupin hypersurface~\cite{Pi3}, and hence 
can be reconstructed from lower dimensional 
proper Dupin hypersurfaces, where the possibilities
for these are already known. The final 
step of the argument uses an unpublished 
result of Pinkall stating that a 
taut hypersurface is analytic.

\medskip

We finish this text with the following lemma which refers to the 
case of a taut hypersurface in $S^5$ with generically
four distinct principal curvatures. 

\begin{lem}
A compact taut hypersurface $M$ in 
$S^5$ with four distinct principal 
curvatures at some point 
has $b_i\geq2$ for $i=1$, $2$, $3$.
\end{lem}

\Pf Let $\xi\in N_p^1(M)$ with $n(\xi)=4$.
Then the multiplicities of the focal points 
in the direction of $\xi$ are $(1,1,1,1)$.
Let $q$ be a point in the normal ray issuing from 
$p$ in the direction of $\xi$ that
comes after the second focal point and before the third one.
Then $L_q$ has a critical point of index 
$2$ at $p$, so $b_2\geq1$. 
Upon a choice of unit normal vector field 
on $M$,  
there is an identification of $N(M)$
with $M\times\R$. 
Let $F:M\times(-\pi,\pi)\subset N(M)\to S^5$
be the restriction of the endpoint map. The focal set of $M$
consists of the critical values of $F$ and consists
of focal varieties of codimension at least~$2$, 
hence its complement $W$ is connected. Let 
$V=F^{-1}(W)$. Then $F|V:V\to W$ is a local
diffeomorphism. For $(x,t)\in V$, let $j(x,t)$ be the
index of $x$ as a critical point of $L_q$ 
for $q=F(x,t)$. Then $j$ is locally constant. 
Decompose $V$ into a disjoint union of nonempty 
open subsets
\[ V = V_0 \cup V_1 \cup V_2 \cup V_3 \cup V_4, \]
where $V_\iota=\{ (x,t)\in V:j(x,t)=\iota\}$. 
If $b_2=1$, 
then $F:V_2\to W$ is a diffeomorphism.
This implies that $V_2$ is connected,
but it is clear that $V_2$ 
must contain points on both sides
of $M\times\{0\}\subset V_0$, 
so this is a contradiction. 
Hence, $b_2\geq2$. 
Similarly, $b_1\geq2$. \EPf 

\section*{Acknowledgement}

We would like to thank Professor D.~L.~Gon\c calves 
and Professor G.~Thorbergsson
for very useful discussions, and specially
the members of the organizing committee of the II Encuentro 
de Geometr\'\i a Diferencial for their outstanding
hospitality. 

\providecommand{\bysame}{\leavevmode\hbox to3em{\hrulefill}\thinspace}
\providecommand{\MR}{\relax\ifhmode\unskip\space\fi MR }
\providecommand{\MRhref}[2]{%
  \href{http://www.ams.org/mathscinet-getitem?mr=#1}{#2}
}
\providecommand{\href}[2]{#2}


\end{document}